\documentstyle[12pt]{article}
\newcommand{\qed} {\hspace {0.1in} \rule {1.5mm} {3.5mm}}

\newtheorem{lemma}{Lemma}
\newtheorem{corollary}[lemma]{Corollary}
\newtheorem{theorem}[lemma]{Theorem}
\newtheorem{claim}[lemma]{Claim}
\def\E{\vec E}
\def\e{\bar e}
\def\d{\delta}
\def\s{\sigma}
\def\R{I\!\!R}
\def\q{\diamondsuit}
\def\t{\star}
\def\dim{{\rm dim}}
\def\supp{{\rm supp}}
\def\<{\langle}
\def\>{\rangle}
\def\proof{\smallskip\noindent{\bf Proof:} }

\title{On quasi-transitive amenable graphs}
\author{{\sc G\'abor Elek} and \sc G\'abor Tardos\cr Mathematical Institute of
the Hungarian Academy of Sciences\cr P.O. Box 127, H-1364 Budapest, Hungary}
\date{}
\begin{document}

\maketitle

\noindent{\bf AMS Subject Classification:} 58G05

\noindent{\bf Keywords:} amenable graphs, harmonic functions

\noindent{\bf Abstract.} We introduce the notion of quasi-transitivity and
prove that there exist no non-constant
harmonic Dirichlet functions on amenable quasi-transitive graphs.

\section{Introduction}
In this note we use the term {\it graph} for simple, connected, undirected
graphs with bounded degree vertices only. We will mostly consider infinite
graphs. We denote the set of vertices of a graph $G$ by $V(G)$, we use
$E(G)$ for the set of edges.
By $\E(G)$ we denote the set
of oriented edges: $\E(G)=\{(x,y)\mid\{x,y\}\in E(G)\}$. We denote the opposite
orientation $(y,x)$ of an oriented edge $e=(x,y)$ by $\e$.
We consider the Hilbert space $l_2(G)$ of the $l_2$
functions $u:\E(G)\to \R$ satisfying $u(\e)=-u(e)$ with the scalar product
$\<u,u'\>=1/2\sum_{e\in\E(G)}u(e)u'(e)$. For simplicity we write $u(x,y)$ for
the value $u((x,y))$. For a function $v:V(G)\to \R$ we define its
{\it differential} $dv:\E(G)\to \R$ by
$dv(x,y)=v(y)-v(x)$. We call a $v:V(G)\to \R$ function a Dirichlet function
if $dv\in l_2(G)$ and denote the set of Dirichlet functions by $D(G)$. Let
$l_2(V(G))$ denote the Hilbert space of the $l_2$ functions $v:V(G)\to \R$
(with the standard scalar product), this is clearly contained in $D(G)$.

Consider the adjoint $d^*$ of the operator $d:l_2(V(G))\to l_2(G)$. We call
a function $u\in l_2(G)$ a {\it flow} if $d^*u=0$. We call a Dirichlet
function $v\in D(G)$ {\it harmonic} if $dv$ is a flow. We denote the set of
harmonic Dirichlet functions by $HD(G)$.

Here $d^*u$ is given by $d^*u(x)=\sum_{\{y,x\}\in E(G)}u(y,x)$ for $u\in
l_2(G)$ and $x\in V(G)$. Thus a $u$ is a flow if and only if $\sum_{\{x,y\}\in
E(G)}u(x,y)=0$ for every $x\in V$. The function $v\in D(G)$ is harmonic if and
only if for every $x\in V$ the value $v(x)$ is the average of the values $v(y)$
with $\{x,y\}\in E(G)$. All constant functions $V(G)\to \R$ are harmonic
Dirichlet functions.

For vertices $x$ and $y$ of a graph $G$ let $\d(x,y)$ denote their
distance in $G$. A {\it wobbling} is a map $f:V(G)\to V(G)$ such that
$\d(x,f(x))$ for $x\in V(G)$ is bounded.
The map $f:V(G)\to V(G')$ is called a {\it quasi-isomorphism} from $G$
to $G'$ if there exits a positive number $k$---the {\it distortion} of
$f$---such that for vertices $x$ and $y$ in $V(G)$ one has
$$\frac1k\d(x,y)-1<\d(f(x),f(y))\le k\d(x,y),$$
and for every vertex $x\in V(G')$ there exists $y\in V(G)$ with
$\d(x,f(y))<k$. A {\it quasi-inverse} of a quasi-isomorphism $f$ from $G$ to
$G'$ is a quasi-isomorphism $g$ from $G'$ to $G$ such that $f\circ g$ and
$g\circ f$ are wobblings.

Note that for a quasi-isomorphism $f$ of distortion $k$ one can take
a quasi-inverse of distortion $2k^2$.

We call a graph $G$ {\it quasi-transitive} if there exist quasi-isomorphisms
$f_{xy}$ from $G$ to $G$ for vertices $x$ and $y$ of $G$ with $f_{xy}(x)=y$
such that these quasi-isomorphisms have bounded distortion.

Any graph quasi-isometric to a vertex-transitive graph is clearly
quasi-transitive. However, the converse is far from being true.
For example, any {\it net} in a non-compact Lie-group or a homogeneous
Riemannian manifold is quasi-transitive. Here the vertices of the net is a
maximal subset of the metric space with minimum distance $1$, the edges connect
vertices of distance at most $3$.

When speaking of subgraphs of a graph we always mean a connected full
subgraph with at least one edge. Let $G_0$ be a subgraph of a graph $G$. By
$\s(G_0)$ we denote the set of vertices of $G_0$ that have neighbors in $G$
outside $G_0$. We call $G$ {\it amenable} if $\inf|\s(G_0)|/|V(G_0)|=0$,
where the infimum is taken for finite subgraphs $G_0$. We call the sequence
$(G_i)$ of finite subgraphs of a graph $G$ a {\it F\o lner
sequence} if $|\s(G_i)|/|V(G_i)|$ tends to $0$ as $k$ tends to infinity.

The study of harmonic Dirichlet functions on graphs goes back to
Cheeger and Gromov \cite{cg}. They proved that there exist no non-constant
harmonic Dirichlet function on the Cayley graph of an amenable group.
For different proofs see also Elek \cite{e} and Paschke \cite{p}.
Later Medolla and Soardi \cite{ms} extended this result to amenable
vertex-transitive graphs. In a recent preprint Benjamini, Lyons, Peres,
and Schramm \cite{blps} gave a probabilistic proof for this result.
The goal of this paper is to extend the result of Medolla and Soardi
to amenable quasi-transitive graphs. Note that some assumption of this kind
is necessary as there exist amenable graphs with non-constant harmonic
Dirichlet functions, see e.g. \cite{s}, Chapter 6.

\section{The result}

We borrow some notations from \cite{blps}.
The {\it support} $\supp(v)$ of a real-valued function $v$ is the subset of
the domain where $v$ is not zero. For a graph $G$ we define $\t(G)$ to be
the closure in $l_2(G)$ of the functions $dv$, where
$v:V(G)\to \R$ has finite support. Let $\q(G)$ be the closure in $l_2(G)$
of the flows with finite support.

Notice that we always have $dl_2(V(G))\subseteq\t(G)$.
The following lemma is well known. We prove it to be self contained.

\begin{lemma}\label1Let $G$ be any graph. We have the following orthogonal
decomposition:
$$l_2(G)=\t(G)+\q(G)+dHD(G).$$
The flows constitute the subspace $\q(G)+dHD(G)$, and $dD(G)=\t(G)+dHD(G)$.
\end{lemma}

\proof First note that if $v\in D(G)$ and $u$ is a flow with finite support
then we have $\<dv,u\>=\<dv_0,u\>=\<v_0,d^*u\>=0$, where $v_0$ is $v$
restricted to the finite support consisting of the endpoints of the
oriented edges in the support
of $u$. Thus $dD(G)$ is orthogonal to $\q(G)$. To see that they are orthogonal
complements consider any function $u$ orthogonal to $\q(G)$. Fix a vertex
$x_0\in V(G)$ and define $v:V(G)\to \R$ by $v(x_n)=\sum_{i=1}^nu(x_{i-1},x_i)$ for
any path $(x_0,\ldots,x_n)$ ($n\ge0$) in $G$. The orthogonality proves that
$v$ is well defined. We have $dv=u$ thus $d\in D(G)$ as required.

Next we claim that the orthogonal complement of $\t(G)$ consists of all the
flows. Indeed a function $u\in l_2(G)$ belongs to this complement if and only if
$\<dv,u\>=\<v,d^*u\>=0$ for every $v:V(G)\to \R$ with a finite support. This
is satisfied if and only if $d^*u=0$.

To finish the proof of the lemma notice that the orthogonal complement of
$\t(G)$ contains $\q(G)$ and the intersection of the orthogonal complements of
$\t(G)$ and $\q(G)$ are flows in $dD(G)$, and by definition, this is $dHD(G)$. 
\qed

Let $G_0$ be a finite subgraph of $G$ and $L$ a closed subspace of $l_2(G)$.
We define
$$\dim_{G_0}(L)={\sum_{e\in\E(G_0)}\<P_Le,e\>\over|\E(G_0)|},$$
where $P_L$ is the orthogonal projection to $L$ and the oriented edge $e$
is identified with the element of $l_2(G)$ mapping $e$ to $1$, $\e$ to $-1$
and all other oriented edges to $0$.

\begin{claim}\label2Let $G_0$ be a finite subgraph of a graph $G$. We have
$\dim_{G_0}(l_2(G))=1$. For orthogonal closed subspaces $L$ and $L'$ of
$l_2(G)$ we have $\dim_{G_0}(L+L')=\dim_{G_0}(L)+\dim_{G_0}(L')$. For closed
subspaces $L\subseteq L'$ one has $0\le\dim_{G_0}(L)\le\dim_{G_0}(L')$. If
the support of all $u\in L$ is contained in $\E(G_0)$ then
$\dim_{G_0}(L)=\dim(L)/|E(G_0)|$.
\end{claim}

\proof
The first statement follows from noting that $P_{l_2(G)}$ is the identity.
The second statement follows from the equality $P_{L+L'}=P_L+P_{L'}$ for
orthogonal subspaces $L$ and $L'$. The non negativity is trivial and implies
the monotonicity. For the last statement note that we can work in $l_2(G_0)$.
Half of the vectors $e\in\E(G_0)$ (taking one of the pairs $e$ and
$\e$) form an orthonormal bases of $l_2(G_0)$, and $|E(G_0)|\dim_{G_0}(L)$
is the trace of $P_L$ restricted to $l_2(G_0)$ computed in this bases.
\qed

\begin{lemma}\label3If $G_0$ is a finite subgraph of the graph $G$ then
$\dim_{G_0}(\t(G)+\q(G))\ge1-|\s(G_0)|/|E(G_0)|$.
\end{lemma}

\proof First note that $dHD(G_0)=0$ since every flow on the finite graph
$G_0$ has finite support, so by Lemma \ref1 $\q(G_0)=\q(G_0)+dHD(G_0)$. So by
Lemma \ref1 we have $l_2(G_0)=\t(G_0)+\q(G_0)$.

We identify functions $u\in l_2(G_0)$ and $v\in D(G_0)$ with
their extension in $l_2(G)$ (respectively $D(G)$) that is zero outside the
original domain. We have to distinguish two differential operators:
$d_{G_0}:D(G_0)\to l_2(G_0)$ is not the restriction of $d:D(G)\to l_2(G)$
as $d(D(G_0))\not\subseteq l_2(G_0)$ unless $G=G_0$.
But using the formula for $d^*$ one sees that $d_{G_0}^*$ is the
restriction of $d^*$. Thus $\q(G_0)\subseteq\q(G)$. Let
$D_1=\{v\in D(G)\mid\supp(v)\subseteq\s(G_0)\}$
and $D_2=\{v\in D(G)\mid\supp(v)\subseteq V(G_0)\setminus\s(G_0)\}$. Clearly,
$D(G_0)=D_1+D_2$ and $\t(G_0)=d_{G_0}D_1+d_{G_0}D_2$. Notice that $d$
and $d_{G_0}$ are identical in $D_2$, thus $d_{G_0}D_2=dD_2\subseteq\t(G)$.

Let $L=dD_2+\q(G_0)$. We have $d_{G_0}D_1+L=\t(G_0)+\q(G_0)=l_2(G_0)$,
thus $\dim(L)\ge\dim(l_2(G_0))-\dim(d_{G_0}D_1)\ge|E(G_0)|-\dim(D_1)\ge
|E(G_0)|-|\s(G_0)|$.

We have $L\subseteq(\t(G)+\q(G))\cap l_2(G_0)$, thus by Claim \ref2
$\dim_{G_0}(\t(G)+\q(G))\ge\dim_{G_0}(L)=\dim(L)/|E(G_0)|\ge
1-|\s(G_0)|/|E(G_0)|$ as claimed.
\qed

\begin{corollary}\label4For a F\o lner sequence $(G_i)$ of
finite subgraphs of a graph $G$ we have $\dim_{G_i}(dHD(G))$ tends to
zero as $i$ tends to infinity.
\end{corollary}

\proof By Lemma \ref1, Claim \ref2, and Lemma \ref3 we have
$0\le\dim_{G_i}(dHD(G))\le|\s(G_i)|/|E(G_i)|$. Here $|\s(G_i)|/|E(G_i)|$ tends
to zero as $(G_i)$ is a F\o lner sequence. Note that for traditional reasons
we used $|\s(G_i)|/|V(G_i)|$ in the definition of F\o lner sequences, but
since $G_i$ is connected and has bounded degree $|V(G_i)|$ and $|E(G_i)|$
are proportional.
\qed

Corollary \ref4 indicates that the harmonic Dirichlet functions on an
amenable graph form a ``small'' subspace. It does not, however, imply that
$dHD(G)=0$, this is false for some amenable graphs.
Our goal is to prove
that if non-constant harmonic functions exist on a quasi-transitive graph
then they form a ``large'' subspace contradicting Corollary \ref4. This is
immediate for transitive graphs: with any harmonic Dirichlet function all its
translates are harmonic. The case of quasi-transitive graphs require more
care. We study next how quasi-isomorphisms act on Dirichlet functions and
on $dHD(G)$.

Let $G_1$ and $G_2$ be graphs. For a map $f:V(G_1)\to V(G_2)$ we define
the function $f^*:v\mapsto f\circ v$ on the functions $v:V(G_2)\to\R$.
Let $A$ be a subset of the vertices of a graph $G$. For $k>0$ we
define the {\it$k$-neighborhood} $C_k(A)$ of $A$ to be $\{x\in
V(G)\mid\exists y\in A:\d(x,y)\le k\}$. We define $\chi_A:\E(G)\to\R$ to
be the characteristic function of the oriented edges $(x,y)\in\E(G)$
with $x\in A$ and $y\in A$.

\begin{lemma}\label5For a quasi-isomorphism $f$ from a graph $G_1$ to a graph
$G_2$ we have $f^*(D(G_2))\subseteq D(G_1)$. Furthermore there is a constant
$c$ depending on the distortion $k$ of $f$ and the maximum degree of $G_1$ such
that $|df^*v|\le c|dv|$ and $|df^*v\cdot\chi_A|\le c|dv\cdot\chi_B|$ for any
$v\in D(G_2)$ and $A\subseteq V(G_1)$ if $B=C_k(f(A))$.
\end{lemma}

\proof Let $f$ be as in the lemma. We fix a path
$f(x)=x_0^e,x_1^e,\ldots,x_{k_e}^e=f(y)$ in $G_2$ of length $0\le k_e\le k$
for each oriented edge $e=(x,y)\in\E(G_1)$. For any $v\in D(G_2)$ we have
$|df^*v|^2=1/2\sum_{(x,y)\in\E(G_1)}(v(f(y))-v(f(x)))^2\le
k/2\sum_{e\in\E(G_1)}\sum_{i=1}^{k_e}(v(x_i^e)-v(x_{i-1}^e))^2$.
The summands in this last expression all appear in the summation
$|dv|^2=1/2\sum_{(x,y)\in\E(G_2)}(v(y)-v(x))^2$ thus for the first statement
we only have to limit the multiplicity of a summand $(v(y)-v(x))^2$ in the
first sum for any $(x,y)\in\E(G_2)$. This is the number of oriented edges in
$\E(G_1)$ such that the corresponding (oriented) path in $G_2$ contains
$(x,y)$. As the endpoints of these edges form a subset in $V(G_1)$ of
maximum distance at most $2k^2$ this multiplicity can be bounded in terms of
$k$ and the maximum degree of $G_1$.

For the last statement notice that if an edge is spanned by a subset
$A\subseteq E(G_2)$ then the corresponding path is within the set
$B=C_k(f(A))$.
\qed

\begin{lemma}\label6For a wobbling $f$ of a graph $G$ we have $f^*v-v\in
l_2(V(G))$ for every $v\in D(G)$.
\end{lemma}

\proof For a vertex $x\in V(G)$ take a path $x=x_0,x_1,\ldots,x_k=f(x)$ in $G$.
We have $((f^*v-v)(x))^2\le k\sum_{i=1}^k(v(x_i)-v(x_{i-1}))^2$. Thus to
bound the $l_2$ norm of $f^*v-v$ in terms of $|dv|$ it is enough to note
that we can choose the paths with bounded length and every edge appears in a
bounded number of paths.
\qed

For a function $f:A\to B$ and $S\subseteq B$ we use $f^{-1}(S)$ to denote
$\{x\in A\mid f(x)\in S\}$.

\begin{lemma}\label7For a graph $G$ we have $dHD(G)\cong D(G)/d^{-1}(\t(G))$.
For a quasi-isomorphism $f$ from $G_1$ to $G_2$\ \ $f^*$ induces an isomorphism
between $D(G_2)/d^{-1}(\t(G_2))$ and $D(G_1)/d^{-1}(\t(G_1))$ thus between
$dHD(G_2)$ and $dHD(G_1)$. For quasi-inverses $f$ and $g$ the functions $f^*$
and $g^*$ induce inverse isomorphisms. 
\end{lemma}

\proof The $dD(G)=\t(G)+dHD(G)$ claim of Lemma \ref1 proves the first
statement.

Let $f$ be as in the lemma and $v\in d^{-1}(\t(G_2))$. Then there exist
functions $v_i\in D(G_2)$ with finite support such that $dv_i$ tend to $dv$
in norm. By Lemma \ref5 $df^*v_i$ tend to $df^*v$, and as the functions
$f^*v_i$ also have finite support $df^*v\in\t(G_1)$. Thus $f^*$ maps
$D(G_2)/d^{-1}(\t(G_2))$ linearly to $D(G_1)/d^{-1}(\t(G_1))$.

To see that this map is an isomorphism take a quasi-inverse $g$ of the
quasi-isomorphism $f$. For any $v\in D(G_2)$ we have by Lemma \ref6 that
$g^*f^*v-v\in l_2(V(G_2))\subseteq d^{-1}(\t(G_2))$. For $v\in D(G_1)$ we
similarly have $f^*g^*v-v\in d^{-1}(\t(G_1))$ thus the maps between
$D(G_2)/d^{-1}(\t(G_2))$ and $D(G_1)/d^{-1}(\t(G_1))$ induced by $f^*$ and
$g^*$ are inverses of each other.
\qed

\begin{lemma}\label8Let the graph $G$ and the positive numbers $k$ and $\Delta$
be given. If $dHD(G)\ne0$ then there exist a finite set $A\subseteq V(G)$ and a
number $\epsilon>0$ with the following property.
For any quasi-isomorphism $f$ of distortion at most $k$ to $G$ from a graph
$G'$ of maximum degree at most $\Delta$ one has a function $w\in HD(G')$ with
$|dw|=1$ and $|dw\cdot\chi_B|>\epsilon$ for $B=f^{-1}(A)$.
\end{lemma}

\proof Choose a non-constant function $v\in HD(G)$. We fix $\epsilon>0$ later
and choose $A=C_{k+1}(A_0)$ with a finite set
$A_0\subseteq V(G)$ such that $|dv\cdot\chi_{V(G)\setminus A_0}|<\epsilon$.

By ``constant'' we mean  a quantity depending on $G$, $v$, $k$ and $\Delta$
but not on $G'$, $f$ or $\epsilon$.
 
We can take a quasi-inverse $g$ of $f$ of distortion bounded by
a constant. By Lemma \ref7 we have the decomposition $f^*v=v'+z$ with some
$v'\in HD(G')$ and $z\in d^{-1}(\t(G'))$. Similarly $g^*v'=v+t$ with some
$t\in d^{-1}(\t(G))$. By Lemma \ref5 we have $0<|dv|\le|dg^*v'|\le c|dv'|$ with
some constant $c$. Thus $|dv'|\ge c_0=|dv|/c>0$. We use Lemma \ref5 for $f^*$
to get $|dz|\le|df^*v|\le c'|dv|=c_1$ with some constants $c'$ and $c_1$.
Now consider $B=f^{-1}(A)\subseteq V(G')$ and $C=(V(G')\setminus B)\cup\s(B)$.
Notice that $\chi_C\ge1-\chi_B$ and $D=C_k(f(C))$ is disjoint from $A_0$.

Consider the orthogonal decompositions $dv'=u_1+u_2$ with $u_1=dv'\cdot\chi_B$
and $u_2=dv'\cdot(1-\chi_B)$, and $dz=s_1+s_2$ with $s_1=dz\cdot\chi_B$ and
$s_2=dz\cdot(1-\chi_B)$. By Lemma \ref5 we have $|u_2+s_2|=
|df^*v\cdot(1-\chi_B)|\le
|df^*v\cdot\chi_C|\le c_2|dv\cdot\chi_D|\le
c_2|dv\cdot\chi_{V(G)\setminus A_0}|< c_2\epsilon$ with a constant $c_2$.
We can write $0=\<dv',dz\>=\<u_1,s_1\>+\<u_2,s_2\>$. Here
$|\<u_1,s_1\>|\le|u_1|\cdot|s_1|\le|u_1|\cdot|dz|\le c_1|u_1|$ and
$\<u_2,s_2\>=\<u_2,u_2+s_2\>-|u_2|^2\le|u_2|\cdot|u_2+s_2|-|u_2|^2<
-|u_2|(|u_2|-c_2\epsilon)$. We have $|u_2|\ge|dv'|-|u_1|\ge c_0-|u_1|$ hence
$c_1|u_1|\ge\<u_1,s_1\>=-\<u_2,s_2\>\ge|u_2|(|u_2|-c_2\epsilon)\ge
(c_0-|u_1|)(c_0-c_2\epsilon-|u_1|)$ if $c_0-c_2\epsilon-|u_1|>0$.

Notice that we can choose a small enough $\epsilon$ depending on $c_0$, $c_1$,
and $c_2$ such that $|u_1|\le c_1\epsilon$ contradicts our last inequality. For
this $\epsilon$ we have $|u_1|>c_1\epsilon$. We take $w=v'/|dv'|\in HD(G')$
and notice that $|dw|=1$ and $|dw\cdot\chi_B|=|u_1|/|dv'|>\epsilon$ as
$|dv'|\le|df^*v|\le c_1$. The choice of $\epsilon$ (and thus of $A$)
depends only on $G$, $v$, $k$ and $\Delta$.
\qed

\begin{theorem}\label9The only harmonic Dirichlet functions of a
quasi-transitive amenable graph are the constant functions.
\end{theorem}

\proof The proof is by contradiction. Let $G$ be a quasi-transitive amenable
graph with $dHD(G)\ne0$. Let $k$ be the bound on the distortion of the
quasi-isomorphisms $f_{xy}$ from $G$ to $G$ mapping $x\in V(G)$ to $y\in V(G)$.
Choose the
finite set $A\subseteq V(G)$ and the number $\epsilon>0$ for $G$, $k$ and the
maximum degree $\Delta$ of $G$ as claimed in Lemma \ref8. We fix a vertex
$y\in A$ and use the statement of Lemma \ref8 to obtain a function $w_x\in
HD(G)$
for each $x\in V(G)$ such that $|dw_x|=1$ and $|dw_x\cdot\chi_{B_x}|>\epsilon$
for $B_x=f_{xy}^{-1}(A)$. Let $a=\max_{y'\in A}\d(y,y')$. By the bound
on the distortion of $f_{xy}$ we have $\d(x,x')<k(a+1)=b$ for any $x\in V(G)$
and $x'\in B_x$.

Let $(G_i)$ be a F\o lner sequence of finite subgraphs of $G$. Let
$S_i=V(G_i)\setminus C_b(\s(G_i))$. Note that as $b$ is constant and the
degree of the vertices is limited, $|C_b(\s(G_i))|$ is proportional to
$|\s(G_i)|$ thus $|S_i|/|V(G_i)|$ tends to $1$ as $i$ tends to infinity.
For $x\in S_i$ we have $B_x\subseteq C_b(\{x\})\subseteq V(G_i)$.
Consider the projection $P$ in $l_2(G)$ to $dHD(G)$ and let $P_x$ the
projection to the one dimensional subspace of $dHD(G)$ generated by $dw_x$.
Recall that we identify an oriented edge $e\in\E(G)$ with the function
in $l_2(G)$ mapping $e$ to $1$, $\e$ to $-1$ and everything else to
zero. We have $\<Pe,e\>\ge\<P_xe,e\>=\<dw_x,e\>^2=(dw_x(e))^2$ for every
oriented edge $e$ and vertex $x$ in $G$. We can write
$\epsilon^2<|dw_x\cdot\chi_{B_x}|^2=1/2\sum_{e\in\E(G)}\chi_{B_x}(e)(dw_x(e))^2
\ge1/2\sum_{e\in\E(G)}\chi_{B_x}(e)\<Pe,e\>$ for any $x\in V(G)$.
By summation we get
$2\epsilon^2|S_i|\le\sum_{e\in\E(G)}\<Pe,e\>\sum_{x\in S_i}\chi_{B_x}(e)$
for any index $i$. Notice that $\sum_{x\in S_i}\chi_{B_x}(e)$ is zero for
oriented edges $e$ outside $\E(G_i)$, while it is bounded by a constant $C$ for
any $e\in\E(G)$, one can take $C$ to be the maximum size of the
$b$-neighborhood of a single vertex. Thus we have
$2\epsilon^2|S_i|/C\le\sum_{e\in\E(G_i)}\<Pe,e\>=|\E(G_i)|\dim_{G_i}(dHD(G))$.
Consequently, $\dim_{G_i}(dHD(G))\ge\epsilon^2/C\cdot|S_i|/|E(G_i)|$.
By Corollary \ref4
the left hand side of this last inequality tends to zero as $i$ tends to
infinity, but as $|S_i|/|V(G_i)|$ tends to $1$ and $|V(G_i)|$ is proportional
to $E(G_i)$ the right hand side does not. The contradiction proves the theorem.
\qed

\end{document}